\documentclass[a4paper,12pt]{amsart}

 \usepackage{graphicx}

\input xy
\xyoption{all}
\usepackage{amssymb}

\usepackage[left=2.5cm,top=2.5cm,right=2.5cm,bottom=2.5cm]{geometry}

\newtheorem{theorem}{Theorem}[section]
\newtheorem{lemma}[theorem]{Lemma}
\newtheorem{proposition}[theorem]{Proposition}
\newtheorem{corollary}[theorem]{Corollary}
\newtheorem*{acknowledgements}{Acknowledgements}
\theoremstyle{definition}
\newtheorem{definition}[theorem]{Definition}

\theoremstyle{remark}
\newtheorem{remark}[theorem]{Remark}
\numberwithin{equation}{section}

\def\R{\mathbb{R}}
\def\N{\mathbb{N}}

\def\T{\mathbb{T}}

\def\F{\mathcal{F}}

\def\A{\mathcal{A}}

\def\H{\mathcal{H}}

\def\-{^{-1}}

\def\K{\mathcal{K}}
\def\V{\mathcal{V}}

\def\C{\mathbb{C}}
\def\Z{\mathbb{Z}}
\def\A{\mathcal{A}}

\begin{document}
\title[Frames generated by actions of countable discrete groups]{Frames generated by actions of countable discrete groups.}

\author{Kjetil R\o ysland}
\footnote{Research supported in part by the Research Council of Norway, project number NFR 154077/420. Some of the final work was also done with support from the project NFR 170620/V30.}
\address{[Kjetil R\o ysland]University of Oslo\\
Department of Biostatistics
Sognsvannsv. 9
PO Box 1122 Blindern
NO-0317 Oslo
Norway\\
Department of Mathematics
PO Box 1053, Blindern
NO-0316 Oslo
Norway
}


\email{roysland@math.uio.no}
\begin{abstract}
We consider dual frames generated by actions of countable discrete groups on a Hilbert space.
Module frames in a class of modules over a group algebra are shown to coincide with a class of ordinary frames 
in a representation of the group. 
 This has applications to shift-invariant spaces and wavelet theory. 
 One of the main findings in this paper is  that whenever a shift-invariant sub space in $L^2( \R^n)$ has compactly supported dual frame generators then it also has compactly supported bi-orthogonal generators.  The crucial part in the proof is a theorem by Swan that states that every finitely generated projective module over the Laurent polynomials in n variables is free. 
\end{abstract}

\maketitle

\section{Introduction.}
In several branches of applied mathematics one seeks suitable decompositions of a Hilbert space $\H$. This could for instance be done using an orthonormal basis, as is often the case in wavelet and Fourier analysis. Moreover, we would also like to be able to construct  our orthonormal basis with some algorithm. An example of such a situation is when $G$ is a countable discrete group with a unitary representation $\nu : G \rightarrow U(\H)$, $\psi \in \H$ is a unit vector and $\{\nu_g \psi\}_{g \in G}$ is an orthonormal basis for $\H$.
However, often one require additional properties such as smoothness, compact support or symmetry properties that are not compatible with orthonormality. If this is the case, then one is quickly led into weakening the assumption of orthonormality. One could, for instance,  try to use frames. A frame in a Hilbert space $\H$ is a family $\{\psi_i  \}_{i \in \N}$ such that there exist $A , B > 0$  such that 
$$
A \| \psi \|^2 \leq \sum_{i \in \N}  |\langle \psi, \psi_i \rangle|^2 \leq B \| \psi \|^2
$$
for every $\psi \in \H$. We say that two frames $\{\psi_i  \}_{i \in \N}$ and 
$\{\tilde \psi_i  \}_{i \in \N}$ are dual if 
$\langle \zeta, \eta \rangle = \sum_{i \in \N} \langle \zeta , \psi_i \rangle \langle \tilde \psi_i , \eta \rangle $ for every $\zeta , \eta \in \H$.
We will consider frames that are generated by a finite family $\psi_1, \dots, \psi_d$ and a unitary representation $\nu$ of 
a countable discrete group $G$, i.e.  such that the vectors  $\{\nu_g \psi_i\}_{g \in G,  1 \leq  i \leq d}$ form a frame in $\H$. 

If we assume an additional decay property  like 
$\langle \psi_i , \nu_g \psi_i \rangle = 0$ for all but finitely many $g \in G$, then it 
turns out that we can answer several questions about such frames using algebraic methods. This is the scope of the present paper. 
The main idea is to consider a dense and $G$-invariant subspace $\Lambda \subset \H$ that consists of elements with a decay property such as the one mentioned above. Now $\Lambda$ is a module over the group algebra $\C(G)$ in a canonical way. The properties of this module will reveal some properties of frames generated  by countable discrete groups.

The idea of using modules in this way to analyze frames is not new. For instance Packer, Rieffel and others have worked on wavelets in the similar setting of projective multiresolution analysis and Hilbert modules, see \cite{PR1}, \cite{Pr2}, \cite{Packer}, \cite{DutRoy1} and \cite{Larsen}. See also \cite{wood}, \cite{hanlarson} and \cite{franklarson}.

The following result is one of the main findings in this paper: Suppose that $\H$ is a shift-invariant subspace in $L^2(\R^n)$, i.e. 
we let $G =\Z^n$ act on $L^2( \R^n)$ by translations and let $\H \subset L^2(\R^n)$ be a closed $G$-invariant subspace. 
 Suppose $\psi_1, \dots, \psi_d$ and $\tilde \psi_1, \dots, \tilde \psi_d$ are compactly supported and 
 $\{\nu_g \psi_i | g \in G , 1 \leq i \leq d\}   $ and $\{\nu_g \tilde \psi_i | g \in G , 1 \leq i \leq d\}   $ form dual frames in  $\H$. 
 In Theorem \ref{main}, we prove that there must exist compactly supported $\phi_1, \dots, \phi_r$ and $\tilde \phi_1, \dots, \tilde \phi_r$ in $\H$ such that 
 $\{\nu_g \psi_i | g \in G , 1 \leq i \leq d\}   $ and $\{\nu_g \tilde \psi_i | g \in G , 1 \leq i \leq d\}   $
  form  bi-orthogonal dual frames for $\H$, i.e. $\langle \nu_{g_1} \phi_i, \nu_{g_2} \tilde \phi_j \rangle = \delta_{g_i,g_j} \delta_{i,j}$.  The crucial part in the proof is a theorem due to Swan, see \cite{Swan}, that states that every finitely generated and projective module over the Laurent polynomials in $n$ variables is also a free module. 

\section{Frames, modules and unitary representations.}

\subsection{Representations of countable discrete groups.}
Let $G$ be a countable discrete group. 
For every $g \in G$, we define a corresponding unitary operator 
$\lambda_g \in B(l^2(G))$ with the equation: 
$$
	(\lambda_g a)_{\tilde g} = a_{g^{-1} \tilde g}.
$$
The operators $\{\lambda_g\}_{g\in G}$ form a  unitary representation of $G$ onto $l^2(G)$.  This is the 
\textit{left regular representation} of $G$. 
We let  $C^*_r(G)$ denote 
the operator norm closure of the 
$*$-algebra generated by the operators $\{\nu_g\}_{g\in G}$. This is called    
the \textit{reduced group} $C^*$-\textit{algebra} 
of $G$. More information about such algebras can be found in  \cite{byexamples}.

Let $f_e$ denote the canonical orthonormal basis element for the trivial element in $G$ and
let $\tau$ denote the state  on $C_r^*(G)$ such that $\tau(a) = \langle a f_e, f_e \rangle$. 
The map  $a \mapsto  \{a_g\}_{g \in G}$ defines an embedding $C_r^*(G) \rightarrow l^2(G)$ 
that translates Parseval's identity into
$
	\tau(aa^*) = \sum_{g \in G} |\tau(a\lambda_g^*)|^2 
$
for every $a \in C^*_r(G)$. Using the polarization identity, we obtain the following useful relation, where the convergence is absolute:  
$$
	\tau(ab^*) = \sum_{g \in G} \tau(a\lambda_g^*)\tau( \lambda_g b^*).
$$
\subsection{Frames}
Let $\H$ be a Hilbert space. We will say that a family  $\{\zeta_i\}_{i \in \N} \subset \H$ is a \textit{frame} for $\H$ if there exist 
$A, B > 0$ such that: 
$$
	A \| \zeta \|^2 \leq \sum_{i \in \N} | \langle \zeta, \zeta_i \rangle| ^2 \leq   B \| \zeta \|^2
$$
for every $\zeta \in \H$. Suppose $\{\tilde \zeta_i \}_{i \in \N}$ is another frame in $\H$. We say that this frame is a \textit{dual frame} of 
 $\{\zeta_i \}_{i \in \N}$ if:  
$$
\zeta = \sum_{i \in \N} \langle \zeta, \zeta_i \rangle \tilde \zeta_i = \sum_{i \in \N} \langle \zeta, \tilde \zeta_i \rangle \zeta_i
$$
for every $\zeta \in \H$. 

Let $\{ e_i\}_{i \in \N}$ denote the canonical orthonormal basis in $l^2( I)$. We define the \textit{analysis operator} $D: \H \rightarrow l^2(\N)$ with the equation:
\begin{align*}
	D\zeta = \sum_{i \in \N} \langle \zeta, \zeta_i \rangle e_i.
\end{align*}
The analysis operator is bounded since $\| D \zeta \|^2 = \sum_{i \in \N} | \langle \zeta , \zeta_i \rangle |^2 \leq B \| \zeta \|^2$. 
Moreover, $D^*D$ is invertible since since 
$ A \|\zeta \|^2 \leq \langle D^*D\zeta,  \zeta\rangle$, see \cite[Proposition 3.2.12]{Pedersen}.   

A standard result states that given a frame $\{\zeta_i\}_{i \in \N}$, then the vectors   
$$
	\{ \tilde \zeta_i = ( D^* D)^{-1} \zeta_i \}_{i \in \N}
$$
form a dual frame. This is  proved in \cite[Proposition 3.2.3]{Daubechies10Lectures}. We say that $\{\tilde \zeta_i\}_{i \in \N}$ is the \textit{canonical dual frame} of $\{\zeta_i\}_{i \in \N}$. 
Finally, if two dual frames $\{\zeta_i\}_{i \in \N}$ and $\{\tilde \zeta_i\}_{i \in \N}$
satisfy $\langle \zeta_i, \tilde \zeta_j \rangle = \delta_{i,j}$, then we will say that they are \textit{bi-orthogonal}.

\subsection{Module frames.}
Let $\A$ be  a unital $*$-algebra and let $E$ be a left $\A$-module. 
A \textit{hermitian form} on $E$ is a map $\langle \cdot, \cdot \rangle_\A : E \times E \rightarrow \A$
such that for every $a \in \A$ and $\zeta_1, \zeta_2 \in E$, then: 
	\begin{enumerate}
			\item $\langle a \zeta_1 + \zeta_2 , \eta \rangle_\A = a \langle \zeta_1, \eta \rangle_\A + \langle \zeta_2, \eta \rangle_\A$,
			\item $\langle \zeta, \eta\rangle_\A^* = \langle \eta, \zeta \rangle_A$.
		\end{enumerate}	
If moreover, $\langle \zeta, \zeta \rangle_\A = 0$ implies that $\zeta = 0$, then we will say that the hermitian form is \textit{non-degenerate}. 
The module is said to be \textit{self dual} with respect to $\langle \cdot , \cdot \rangle_\A$ if for every 
$\omega \in \text{Hom}_\A ( E, \A)$, there exists an $\eta \in E$ such that $\omega = \langle \cdot, \eta \rangle_\A$.

Suppose $\langle \cdot, \cdot \rangle_\A$  is a non-degenerate hermitian form on $E$. 
Two families 
	$\zeta_1, \dots, \zeta_d \in E$ and $\tilde \zeta_1  \dots, \tilde \zeta_d \in E$ form \textit{dual module frames} for $E$ if
  $$
		\langle \zeta , \eta \rangle_\A = \sum_{i = 1}^d \langle \zeta, \zeta_i \rangle_\A	\langle \tilde \zeta_i , \eta \rangle_\A
	$$
	for every $\zeta, \eta  \in E$.  
	Using the properties of  $\langle \cdot , \cdot \rangle_\A$, we can easily check that  $\zeta_1, \dots, \zeta_d $ and $\tilde \zeta_1  \dots, \tilde \zeta_d $ form dual module frames if and only if 
		$$
		\sum_i \langle \zeta, \zeta_i \rangle_\A \tilde \zeta_i = \sum_i \langle \zeta, \tilde \zeta_i \rangle_\A  \zeta_i = \zeta
	$$
for every $\zeta \in E$. 
If in addition $\langle \zeta_i , \tilde \zeta_j \rangle_\A = \delta_{i,j}$ for every $1 \leq i, j \leq d$, 
then we will say that the dual frames form 
\textit{dual module bases}. 

\subsection{Projective modules.}
Recall that a finitely generated  $\A$-module is said to be projective if there exists another $\A$-module $F$ such that 
$E \oplus F$ is free,  i.e. $E \oplus F$ is isomorphic to $\A^m$ for an $m \in \N$.
An $\A$-module is projective  if and only if for every $\A$-module $M$ and surjective $T \in \text{Hom}_\A ( M, E)$ there exists an 
$S \in \text{Hom}_\A ( E, M)$ such that $ST = id_E$, see \cite{Lang}.

\begin{lemma} \label{moduleframelemma}
	A finitely generated $\A$-module $E$ with a non-degenerate hermitian form has dual module frames if and only if it is projective and self dual. 
	Moreover, $E$ has dual module bases if and only if it is free and self dual.
	\begin{proof}
		First we assume that $E$ is self dual. 
		Let $F$ be  another $\A$-module with a hermitian form 
		and let 
  	$S \in \text{Hom}_\A ( E, F)$. Now, $\langle S \cdot, \eta \rangle_\A \in  \text{Hom}_\A ( E, \A)$ for every $\eta \in F$ and there 	
		exists an element $S^*\eta \in E$ such that $\langle S \zeta, \eta \rangle_\A = \langle \zeta, S^* \eta \rangle_\A$ for every $\zeta \in E$. 
		We can easily check that the map $\eta \mapsto S^*\eta$ is $\A$-linear. 
		This implies that a finitely generated module $E$ with a non-degenerate and hermitian form is self dual if and only if 
		for every finitely generated $\A$-module $F$ with a hermitian form then 
		every  $S \in \text{Hom}_\A ( E, F)$ is adjointable, i.e. there exists  an $S^* \in  \text{Hom}_\A ( F, E)$ such that 
		$\langle S \zeta, \eta \rangle_\A = \langle \zeta, S^* \eta \rangle_\A$ for every $\zeta \in E$ and $\eta \in F$.
	
		Suppose $E$ is finitely generated,  projective and self dual. We let $\zeta_1, \dots, \zeta_n$ be a generating set for $E$ and 
		let $S$ denote the $\A$-linear map $\A^d \rightarrow E$ such that $Se_i = \zeta_i$ for $1 \leq i \leq d$. Since $E$ is projective, there exists 
		a $T \in \text{Hom}_\A( E, \A^d)$ such that $ST = 1_E$.
		Now, we equip 
		$\A^d$ with the hermitian form $\langle a, b \rangle_\A = \sum_{i = 1}^d a_i b_i^*$. The map $T$ is adjointable since $E$ is self dual. 
		Let $\tilde \zeta_i = T^* e_i$ for $1 \leq i \leq d$.  We see that whenever $\zeta \in E$, then
		\begin{align*}
			\sum_{i = 1}^d \langle \zeta, \tilde \zeta_i \rangle_\A \zeta_i & = \sum_{i = 1}^d \langle \zeta, T^* \zeta_i \rangle_\A S e_i
			= S \sum_{i = 1}^d \langle T \zeta,  \zeta_i \rangle_\A  e_i = ST \zeta = \zeta,
		\end{align*}
		so $\zeta_1, \dots, \zeta_d $ and $\tilde \zeta_1  \dots, \tilde \zeta_d $ form dual module frames for $E$. Suppose moreover, that $E$ is free 
		and assume that $\zeta_1, \dots, \zeta_d$ form an $\A$-basis for $E$. Now $S$ is invertible, so if we define 
		$
			\tilde \zeta_i = (S^{-1})^* e_i, 
		$	
		then we see that 
		\begin{align*}
			\langle \zeta_i , \tilde \zeta_j \rangle_\A = \langle \zeta_i, (S^{-1})^* e_i \rangle_\A = \langle S^{-1} \zeta_i , e_j \rangle_\A
			= \langle e_i, e_j \rangle_\A = \delta_{i,j}. 
		\end{align*}

	For the converse statement,  suppose $\zeta_1, \dots, \zeta_d $ and $\tilde \zeta_1  \dots, \tilde \zeta_d $ form dual module frames for $E$ and let 
		$\omega \in \text{Hom}_\A(E, \A)$. If $\zeta \in E$, then 
		\begin{align*}
			\omega(\zeta) = \omega( \sum_{i = 1}^d \langle \zeta , \tilde \zeta_i \rangle_\A \zeta_i ) = \langle \zeta, \sum_{i = 1}^d 
			\omega(\zeta_i)^* \tilde \zeta_i \rangle_\A,
		\end{align*}
	so $E$ is self dual. Now, we define  $S, \tilde S \in \text{Hom}_\A( \A^d, E)$  such that $S e_i = \zeta_i$ and $\tilde S e_i = \tilde \zeta_i$ for $ 1 \leq i \leq d$. A short computation  shows that $S$ and $\tilde S$ are adjointable and 
	\begin{align*}
	& S^* \zeta  = \sum_{i = 1}^d \langle \zeta, e_i \rangle_\A \zeta_i 
	& \tilde S^* \zeta  = \sum_{i = 1}^d \langle \zeta, e_i \rangle_\A \tilde \zeta_i. 
	\end{align*}	 
 Moreover, we can easily check that $S \tilde S^* = \tilde S S^* = 1_E$. This implies that 
 $$\zeta \mapsto S \zeta \oplus ( 1 - \tilde S^* S)  \zeta $$ defines an $\A$-linear isomorphism from $\A^d$ to $E \oplus \ker \tilde S^* S$, so  $E$ is projective.  
		
	\end{proof}
\end{lemma}

\section{The module of compactly supported vectors in $\H$.}
 We let $\A$ denote the $*$-algebra of finite linear combinations of elements in 
$\{\lambda_g\}_{g \in G} \subset C_r^*(G)$. 
Moreover, let $\H$ be a Hilbert space with a unitary representation $\nu: G \rightarrow U( \H)$. 

The action  $\nu$ defines a $*$-homomorphism $\pi_0 : \A \rightarrow B(\H)$ such that 
$\pi_0( \lambda_g)  = \nu_g$. However,  we can not necessarily extend this to $*$-homomorphism $C^*_r(G)\rightarrow B( \H)$, unless $G$ is amenable, in which case it is well known that $C_r(G)$ coincides with the full group $C^*$-algebra on $G$. 
We will occasionally apply the operator norm inherited from $C_r^*(G)$ on elements in $\A$. This norm will be denoted by $\| \cdot \|$.
Moreover, the restriction of the cone of positive elements in $C_r^*(G)$ defines a cone of positive elements in $\A$.

We let  $\Lambda \subset \H$ be a subspace such that:
\begin{enumerate}
\item $\pi_0(a) \Lambda \subset \Lambda$ for every $a \in \A$, 
\item If $\zeta, \eta \in \Lambda$, then $\langle \zeta, \nu_g \eta \rangle = 0$ for all but finitely many $g \in G$.  \label{decaymodule}
\end{enumerate}
The first condition means  that $\Lambda$ is an $\A$-module when equipped 
with the product $a, \zeta \mapsto \pi_0(a) \zeta$. 
The motivating example for this is when $\Lambda$ is a subspace of $L^2(\R^n)$ that is formed by functions of compact support and 
$G = \Z^n \subset \R^n $ acts on $L^2(\R^n)$ by translations. Note, that most computations and definitions work out nicely if we weaken the decay assumption, i.e. instead of assuming (\ref{decaymodule}), we  assume that $\{\langle \zeta, \nu_g \eta \rangle\}_{g \in G} \in l^1( G)$.  
 
Now, we define a hermitian form $\langle \cdot ,  \cdot \rangle_\A$ on the $\A$-module $\Lambda$ with the following equation:
\begin{align}
	\langle \zeta, \eta \rangle_\A = \sum_{g \in G} \langle \zeta, \nu_g \eta \rangle \lambda_g.
\end{align}
We list some of its properties in the  following proposition:  
\begin{proposition} \label{inner product}
Let  $\zeta_1, \zeta_2, \zeta , \eta \in \Lambda$ and $a\in \A$. 	
		The map $\langle \cdot ,  \cdot \rangle_\A$ satisfies the following properties:  
		\begin{enumerate}
			\item $\langle \zeta, \zeta \rangle_\A \geq 0 $ and $\langle \zeta, \zeta \rangle_\A = 0$ if and only if  $\zeta = 0$.
			\item $\langle a \zeta_1 + \zeta_2 , \eta \rangle_\A = a \langle \zeta_1, \eta \rangle_\A + \langle \zeta_2, \eta \rangle_A$. 
			\item $\langle \zeta, \eta\rangle_\A^* = \langle \eta, \zeta \rangle_A$.
			\item \label{Cauchy-Schwartz}
			$\langle \zeta, \eta \rangle_A \langle \eta , \zeta \rangle_\A \leq \| \langle \zeta, \zeta \rangle_\A \| \langle \eta , \eta \rangle_\A$.
		\end{enumerate}	
		\begin{proof}
					We let $a, b \in l^2(G)$ and compute with respect to the standard inner product on $l^2(G)$: 
			\begin{align*}
					\langle \langle \zeta, \eta \rangle_\A a, b \rangle & = 
				\sum_{g \in G}\sum_{f \in G} \langle \zeta, \nu_f \eta \rangle a_{f^{-1} g} \overline a_g 
				 = \sum_{g \in G}\sum_{f \in G} \langle  \zeta, \nu_{f^-1} \eta \rangle a_{g} \overline a_{f^-1}g \\
				& = \langle a, \langle \eta, \zeta \rangle_\A b \rangle. 
			\end{align*}
			This shows that $\langle \zeta, \eta \rangle^*_\A = \langle \eta , \zeta \rangle_\A$ for every $\zeta, \eta \in \Lambda$.
			
			To see that $\langle \zeta, \zeta \rangle _\A \geq 0$, we pick  $a \in l^2(G)$ with 
			finite support and consider the following computation:  
			\begin{align*}
			\langle \langle \zeta, \zeta \rangle_\A a, a \rangle & = 
			\sum_{g \in G} \sum_{f \in G} \langle \zeta , \nu_f \zeta \rangle a_{f^{-1}g} \overline a_g 
			 = \sum_{f \in G} \sum_{g \in G} \langle \zeta, \nu_{g^{-1} f} \zeta \rangle a_g \overline a_f \\
				& = \langle \sum_{g \in G} a_g \nu_g \zeta, \sum_{f \in G} a_f \nu_f \zeta \rangle  \geq 0.  
			\end{align*}
			Now,  let $a \in l^2(G)$ be arbitrary and let $\{a_n\}_{n \in \N} \subset l^2(G)$ be a sequence 
			of elements with finite supports such that $\lim_n \| a_n -a \| = 0$. We see that 
			$\langle \langle \zeta, \zeta \rangle_\A a, a\rangle = \lim_n \langle \langle \zeta, \zeta \rangle_\A a_n, a_n\rangle \geq 0$, 
			so $\langle \zeta, \zeta \rangle_\A\geq 0$ for every $\zeta \in \Lambda$.   Moreover, if $\langle \zeta , \zeta \rangle_\A = 0$, then 
			$\langle \zeta, \nu_g \zeta \rangle = 0$ for every $g \in G$ since $\{\lambda_g\}_g$ yields an orthonormal basis for $l^2(G)$. 
			
			Whenever $a\in \A$, we have the following relation:
			\begin{align*}
				\langle a \zeta, \eta \rangle_\A & = \sum_{g \in G} \sum_{f \in G} a_f \langle \nu_g \zeta, \nu_f\eta \rangle \lambda_f
					= \sum_{g \in G} \sum_{f \in G} a_g \lambda_g \langle \zeta, \nu_f \eta \rangle \lambda_f = a \langle \zeta, \eta \rangle_\A.
			\end{align*} 
			Now, its easy to see that $\langle a \zeta_1 + \zeta_2 , \eta \rangle_\A = a\langle \zeta_1, \eta \rangle_\A + \langle \zeta_2, \eta \rangle_\A$
			whenever $\zeta_1, \zeta_2 , \eta \in \Lambda$. 
			
			Finally, we will prove the Cauchy-Schwartz like inequality (\ref{Cauchy-Schwartz}) as follows:
			\begin{align*}
				0 & \leq 
				\langle \langle \zeta, \eta \rangle_\A \| \langle \eta, \eta \rangle_\A \|^{-1} \eta - \zeta , 
				\langle \zeta, \eta \rangle_\A \| \langle \eta, \eta \rangle_\A \|^{-1} \eta - \zeta \rangle_\A \\
					& = \langle \zeta, \eta\rangle_\A \langle \eta, \eta \rangle_\A \langle \eta, \zeta \rangle_\A 
						\| \langle \eta, \eta \rangle_\A \|^{-2} - 2 \langle \zeta , \eta \rangle_\A \langle \eta, \zeta \rangle_\A \| \langle \eta, \eta \rangle_\A \|^{-1} + \langle \zeta, \zeta \rangle_\A. 
			\end{align*} 
			If $a,b$ are elements in a $C^*$-algebra and $b$ is positive, then  $a^*ba \leq \|b\| a^*a$, see \cite[2.2.2]{BratteliRobinson}. Since $\langle \zeta , \zeta \rangle_\A$ is positive in $C^*_r( G)$, we obtain
			the inequality $$\langle \zeta, \eta \rangle_\A \langle \eta, \eta \rangle_\A \langle \eta , \zeta \rangle_\A 
				\leq  \| \langle \eta, \eta \rangle_\A \|^{-1} \langle \zeta, \eta \rangle_\A \langle \eta , \zeta \rangle_\A.$$ Using this inequality, we obtain the following: 
			\begin{align*}	
				0 \leq - \langle \zeta, \eta \rangle_\A \langle \eta, \zeta \rangle_\A \| \langle \eta , \eta \rangle_\A \|^{-1} + \langle \zeta, \zeta \rangle_\A.
			\end{align*}
			We see immediately  that this is equivalent to (\ref{Cauchy-Schwartz}).
		\end{proof}
\end{proposition}
 
\begin{remark}
	The module $\Lambda$ is closely related to a Hilbert module introduced by  Packer and Rieffel in \cite{PR1}. 
 They considered the space: 
	$$\Xi  = \{ f \in C_b ( \R^n) | \sum_{ g \in \Z^n   } f ( \cdot - g ) \overline f (\cdot - g) \text{
	converges uniformly.}\}.$$
	This space $\Xi$ is a $C(\T^n)$-module in a canonical way. It is also a $C(\T^n)$-Hilbert module, see \cite{Lance}, when equipped with 
	the $C(\T^n)$-valued hermitian form:  
	 $$\langle f_1, f_2 \rangle'  = \sum_{ g \in \Z^n  } f_1 ( \cdot - g ) \overline f_2 (\cdot - g).$$  
	They introduced this Hilbert module in order to define  
 	\textit{projective multiresolution analysis}, see also \cite{DutRoy1} and \cite{Larsen}.
	We will not deal with projective resolution analysis, but only make a few comments on the relation between $\Lambda$ and 
	the Hilbert module $\Xi$.
	
	Suppose,  as in the previous section,  that $G = \Z^n$ and assume that the action $\nu$ on $L^2( \R^n)$ is given as follows:
	$
		(\nu_g \zeta)(x) = \zeta( x - g), g \in G.
	$
	We let $\Lambda$ denote the set of compactly supported functions in $L^2(\R^n)$.  
	The Fourier transform relates the $\A$-module $\Lambda$ to the Hilbert module $\Xi$ in the following way:  
	First, we note that if $\zeta \in \Lambda \subset L^1(\R^n)$, then $\widehat \zeta$ is  continuous. 
	Let $G^\bot = \{y \in \R^n | \langle x,y \rangle \in \Z \text{ for every } x \in G\}$ and 
	let $X = \R^n / G^\bot$. Moreover, let $\mu$ be the Haar measure on $X \simeq \T^n$ with the normalization such that 
	$
	 \int_{\R^n} f dx = \int_X \sum_{g \in G^\bot} f( x - g) d\mu([x])	
	$ for every $f \in C_c( \R^n )$. Now if $\zeta_1, \zeta_2, \eta  \in \Lambda$,  then  
	\begin{align*}
	& \widehat{   (\langle \zeta_1, \zeta_2 \rangle_\A \eta )}( \sigma)  =  
	\sum_{g \in G} \langle \zeta_1, \nu_g \zeta_2 \rangle \widehat{\nu_g \eta}( \sigma)    =      \sum_{g \in G} \int_{\R^n} \hat \zeta_1(y) \overline{ \hat \zeta_2} ( y ) \overline{(g, y)} dy  (g,\sigma )  \hat \eta (\sigma) \\
	 = & (\F^{-1} \F \langle \zeta_1, \zeta_2 \rangle')(\sigma)\hat \eta(\sigma) 
	 =   \langle \hat \zeta_1, \hat \zeta_2 \rangle' ( \sigma) \hat \eta ( \sigma). 
\end{align*}	
	This means that the Fourier transform takes $\Lambda$ to a subspace in $\Xi$ and  
	takes the hermitian form $\langle \cdot , \cdot \rangle_\A$
	to $\langle \cdot , \cdot \rangle'$. 
\end{remark}
	
\begin{remark}	
A countable group $G$ is called \textit{amenable} if there exists a left translation invariant state on $l^\infty (G)$. 
Finite groups and Abelian groups are examples of amenable groups.  This is a standard result that is proved in \cite[VII.2]{byexamples}. Moreover, if 
 $G$ is  an amenable discrete group, $\nu$ is a unitary representation of  $G$ on $\H$ and $\mathcal{B}$ is the $C^*$-algebra generated by the representation $\nu$, then there exists a  $*$-homomorphism $\pi : C^*_r(G) \rightarrow \mathcal{B}$ such that $\pi(\lambda_g) = \nu_g$ for every $g \in G$.
This result is also standard and can be found in 
 \cite[Theorem V.II.2.8]{byexamples}.  
 We note that if such a $\pi$ exists, then it is unique. This follows from extension by continuity since 
$\A$ is a dense $*$-sub algebra in $C_r^*(G)$. 

If $G$ is amenable, then we can always complete $\Lambda$ into a $C_r^*(G)$-Hilbert module $\Lambda_0 \subset \H$.  
Let $\tilde \Lambda$ denote the vector space $\pi(C_r^* ( G)) \Lambda$.
The space $\tilde \Lambda$ has a canonical $C_r^*(G)$-module structure and  the hermitian form $\langle \cdot, \cdot \rangle_\A$ extends uniquely to a $C_r^*(G)$-valued hermitian form $\langle \cdot, \cdot \rangle_\A$ on $\tilde \Lambda$. 
The map $\zeta \mapsto \| \langle \zeta, \zeta \rangle _{C_r^* (G)}\|^{1/2}_{\text{op}}$ defines a norm on $\tilde \Lambda$, 
so the usual completion procedure of a normed vector space defines a left module over $C_r^*(G)$. 
If a sequence $\{\zeta_n\}_n \subset \tilde \Lambda$ is a Cauchy sequence with respect to this norm, then 
it is also a Cauchy sequence with respect to the norm on $\H$. Since $\H$ is complete, we see that 
the completion of $\tilde \Lambda$ is contained in $\H$. We let $\Lambda_0\subset \H$ denote the subspace formed by this completion.   
The hermitian form $\langle \cdot , \cdot \rangle_\A$ now extends to a non-degenerate $C_r^*(G)$-valued hermitian form 
$\langle \cdot , \cdot \rangle_{C_r^* (G)} : \Lambda_0 \times \Lambda_0 \rightarrow C_r^* (G)$ that makes $\Lambda_0$ into a Hilbert module, i.e. 
$\langle \zeta  , \zeta \rangle_{C_r^* (G)} \geq 0$ for every $\zeta \in \Lambda_0$ and 
$\langle \zeta, \eta \rangle_{C_r^*(G)} \langle \eta , \zeta \rangle_{C_r^*(G)} \leq \| \langle \zeta, \zeta \rangle_{C_r^*(G)} \| \langle \eta , \eta \rangle_{C_r^*(G)}$ for every $\zeta, \eta \in \Lambda_0$.
		
\end{remark}

\subsection{$G$-frames.}


We will say that a finite family of vectors $\zeta_1, \dots ,\zeta_d \in \H$ form a $G$-\textit{frame} in $\H$ 
if  the vectors $\{\nu_g \zeta_i | g \in G , 1 \leq i \leq d\}$ form a frame in $\H$. The following theorem says that we can tell if two families
$\{\zeta_i\}_{ 1 \leq i \leq d}$ and $\{\tilde \zeta_i\}_{ 1 \leq i \leq d}$
form dual $G$-frames for a closed subspace in $\H$ by looking at a corresponding sub module in $\Lambda$. A similar result can be found in 
\cite{PR1}.

\begin{theorem} \label{frametheorem}Suppose  that $E \subset \Lambda$ is an $\A$-submodule and let 
		$\zeta_1, \dots, \zeta_d , \tilde \zeta_1  \dots, \tilde \zeta_d \in E$. The following conditions are equivalent:
		\begin{enumerate}
			\item The families $\{\zeta_i\}_{ 1 \leq i \leq d}$ and $\{\tilde \zeta_i\}_{ 1 \leq i \leq d}$ form dual $G$-frames for $\overline E$.
			\item The families $\{\zeta_i\}_{ 1 \leq i \leq d}$ and $\{\tilde \zeta_i\}_{ 1 \leq i \leq d}$ form dual module frames for $E$.
		\end{enumerate}
		\begin{proof}

		Suppose $\zeta\in E$	and that the families 
	 	$\{\zeta_i\}_{ 1 \leq i \leq d}$ and $\{\zeta_i\}_{ 1 \leq i \leq d}$ form dual $G$-frames for the 
			closure of $E$ in $\H$. We see that 
				\begin{align*}	
					\zeta = \sum_{g,i} \langle \zeta, \nu_g \zeta_i \rangle \nu_g \tilde \zeta_i = \sum_i \langle \zeta, \zeta_i \rangle_\A \tilde \zeta_i,
				\end{align*}
		so $\{\zeta_i\}_{ 1 \leq i \leq d}$ and $\{\zeta_i\}_{ 1 \leq i \leq d}$ form dual module frames for $E$.

		If $\zeta, \eta \in \Lambda$, then 
		\begin{align} 
			\sum_{g \in G} \langle  \zeta, \nu_g \eta \rangle \langle \nu_g \eta, \zeta \rangle  & = \sum_{g \in G }\tau( \langle \zeta, \eta \rangle_\A \lambda_g^*)
			    \tau( \lambda_g \langle \eta, \zeta \rangle_\A) 
			    	= \tau( \langle \zeta, \eta \rangle_\A \langle \eta, \zeta \rangle_\A) \\   & \leq \tau( \langle \zeta , \zeta \rangle_\A) \| \langle \eta , \eta \rangle_\A \| = \| \zeta \|^2 \| \langle \eta , \eta \rangle_\A \|. \label{ineq1}
		\end{align}
		This implies that $\{\langle \zeta, \nu_g \eta \rangle\}_{g \in G}$ is contained in $l^2(G)$. Moreover, if 
		$\{\zeta_n\}_n$ is a sequence in $\H$ that converges to $\zeta \in \H$, then  
		$l^2-\lim_n \{\langle \zeta_n, \nu_g\eta \rangle\}_{g \in G} = \{\langle \zeta, \nu_g \eta \rangle\}_{g \in G}$.
		
		Suppose that $\{\zeta_i\}_{ 1 \leq i \leq d}$ and $\{\tilde \zeta_i\}_{ 1 \leq i \leq d}$ form dual module frames for $E$.
		If  $\zeta, \eta $ are  contained in the closure of $E$, let $\{\zeta^{(n)}\}_{n \in \N}$  and $\{\eta^{(n)}\}_{n \in \N}$ be  sequences in $E$ that converge to $\zeta$ and $\eta$.  
		We see that 
		\begin{align*}
			\langle \zeta, \eta \rangle & = \lim_n \langle \zeta^{(n)} , \eta^{(n)} \rangle = \lim_n \tau(\langle \zeta^{(n)}, \eta^{(n)} \rangle_\A) = 
			 \lim_n \tau( \sum_i \langle \zeta^{(n)}, \zeta_i \rangle_\A \langle \tilde \zeta_i, \eta^{(n)} \rangle_\A     ) \\
			 	& =  \lim_n \sum_{i,g} \tau( \langle \zeta^{(n)}, \zeta_i \rangle_\A \lambda_g^*) \tau( \lambda_g \langle \tilde \zeta_i, 
			 	\eta^{(n)} \rangle_\A) 
			 	= \lim_n \sum_{g, i} \langle \zeta^{(n)}, \nu_g \zeta_i \rangle \langle \nu_g \tilde \zeta_i , \eta^{(n)} \rangle \\
			 	& =  \sum_{g, i} \langle \zeta, \nu_g \zeta_i \rangle \langle \nu_g \tilde \zeta_i , \eta \rangle.  
		\end{align*}
		To prove that $\zeta_1, \dots, \zeta_d$ form a $G$-frame for  the closure of $E$ in $\H$ we must show that there exist
		$A,B > 0$ such that 
		\begin{align} \label{ineq2}
		A\| \zeta \|^2 \leq \sum_{g,i} |\langle \zeta, \nu_g \zeta_i \rangle|^2 \leq B \|\zeta \|^2
		\end{align}
		for every $\zeta$  in the closure of $E$.
		The second inequality of (\ref{ineq2}) is satisfied if we set
		$B = \sum_{i} \| \langle \zeta_i,\zeta_i \rangle_\A \|$.
		Moreover, we have that 
			\begin{align*}
				 \| \zeta \|^2 & = \sup_{ \| \eta \| \leq 1} | \langle \zeta, \eta \rangle |^2 = 
				 \sup_{ \| \eta \| \leq 1} | \sum_{g,i} \langle \zeta, \nu_g \zeta_i \rangle \langle \nu_g \tilde \zeta_i , \eta \rangle |^2 \\
				 & \leq \sup_{ \| \eta \| \leq 1} \sum_{g,i} |\langle \zeta, \nu_g \zeta_i \rangle|^2 \sum_{g,i} |\langle \eta, \nu_g \tilde \zeta_i \rangle|^2
				 \leq \sup_{ \| \eta \| \leq 1}
				 	   \sum_{g,i} |\langle \zeta, \nu_g \zeta_i \rangle|^2  \| \eta \|^2 \sum_{i} \| \tilde \langle \zeta_i, \tilde \zeta_i \rangle_\A \| \\
				 & \leq \sum_{g,i} |\langle \zeta, \nu_g \zeta_i \rangle|^2 \sum_{i} \| \langle \tilde  \zeta_i, \tilde \zeta_i \rangle_\A \|.  
			\end{align*}
		So, if we let $A = (\sum_{i} \| \langle\tilde  \zeta_i, \tilde \zeta_i \rangle_\A \|)^{-1}$, then 
		$
			A \| \zeta \|^2 \leq \sum_{g,i} |\langle \zeta, \nu_g \zeta_i \rangle|^2
		$
		for every $\zeta$ in the closure of $E$ in $\H$. Finally, we note that this argument can also be done for the family 
		$\tilde \zeta_1, \dots, \tilde \zeta_d$, so the families 
		$\{\zeta_i\}_{ 1 \leq i \leq d}$ and $\{\zeta_i\}_{ 1 \leq i \leq d}$ form dual $G$-frames for the 
			closure of $E$ in $\H$.
		\end{proof}
\end{theorem}

\begin{corollary}\label{biortcor}
	 Suppose  that $E \subset \Lambda$ is an $\A$-submodule and let 
		$\zeta_1, \dots, \zeta_d , \tilde \zeta_1  \dots, \tilde \zeta_d \in E$. The following conditions are equivalent:
		\begin{enumerate}
			\item The families $\{\zeta_i\}_{ 1 \leq i \leq d}$ and $\{\zeta_i\}_{ 1 \leq i \leq d}$ form bi-orthogonal $G$-frames for $\overline E$.
			\item The families $\{\zeta_i\}_{ 1 \leq i \leq d}$ and $\{\zeta_i\}_{ 1 \leq i \leq d}$ form dual module bases for $E$.
		\end{enumerate}
	\begin{proof}
		This follows directly from the previous theorem and the definition of $\langle \cdot , \cdot \rangle_\A$.
	\end{proof}
\end{corollary}

\section{Shift-invariant subspaces in $L^2( \R^n)$.}
In this section we will apply our results to shift-invariant subspaces in $L^2( \R^n)$.
We let $G= \Z^n$ and consider  the representation $\nu :G \rightarrow \mathcal{U}(L^2( \R^n))$
such that  
	$$
		(\nu_g \zeta)(x) = \zeta( x - g), g \in G.
	$$
A closed subspace  $\mathcal{V} \subset L^2(\R^n)$ is said to be shift-invariant if 
$\nu_g \mathcal{V} \subset \mathcal{V}$ for every $g \in G$. Note that since $G = \Z^n$, the algebra  $\A$ is isomorphic to the $*$-algebra of  Laurent polynomials  in $n$-variables.

		We let $\Lambda$ denote the set of functions in $L^2(\R^n)$ of compact support. 
		If $\zeta_1, \dots , \zeta_d \in \Lambda$ form a $G$-frame for a closed subspace 
		$\mathcal{V} \subset L^2( \R^n)$,  
		then we will see that the structure of $\V$ is closely related to the structure of the $\A$-module 
		$\A \zeta_1 + \dots + \A \zeta_d$. 
		
		A frame $\{\zeta_i\}_{i \in \N}$ in a Hilbert space $\H$ is said to be a \textit{Parseval frame} for $\H$ if 
		$$ \langle \zeta, \eta \rangle  = \sum_{i \in \N} \langle \zeta, \zeta_i \rangle \langle \zeta_i , \eta \rangle$$
		for every $\zeta, \eta  \in \H$.  
			
		\begin{theorem}
			If $\{\nu_g \zeta\}_{g \in G} \subset  \mathcal {V} \subset L^2(\R^n)$ is a Parseval frame for $\mathcal{V}$ and $\zeta$ has compact support, then 
			$\{\nu_g \zeta\}_{g \in G}$ is an orthonormal basis for $\mathcal{V}$. 
	\begin{proof}
	Let $E$ denote
	 the $\A$-module that consists of the compactly supported elements in $\mathcal{V}$.
	 	By Theorem \ref{frametheorem},  $\zeta$ forms a module frame for  $E$ such that 
	$$\langle \eta_1, \eta_2 \rangle_\A   = \langle \eta_1, \zeta \rangle_\A \langle \zeta, \eta_2 \rangle_\A$$
	for every $\eta_1, \eta_2 \in E$. This implies that 
	$\langle \zeta, \zeta   \rangle_\A^2 = \langle \zeta, \zeta \rangle_\A$. Moreover, $\langle \zeta, \zeta \rangle_\A$ is self-adjoint and non-zero. 
	The Fourier transform yields a $*$-isomorphism between $C_r^*( \Z^n)$ and $C(\T^n)$, see \cite[Proposition VII.1.1]{byexamples}. The only self adjoint idempotents in this algebra are $1$ and $0$, so 
	 $\langle \zeta, \zeta \rangle_\A = 1$. If we apply Corollary  \ref{biortcor}, we see  that 
	$\langle \zeta, \nu_g \zeta \rangle = \delta_{g,0}$. 
		\end{proof}
	\end{theorem}



A fundamental fact about polynomial rings in $n$ variables over a field is that every 
finitely generated and projective module over such a ring is free. This is known as Serre's conjecture and was an open problem for several years until it was proved independently by Quillen and Suslin in 1976. Richard G. Swan proved the analogues result for Laurent Polynomials in $n$ variables in \cite{Swan}, see also \cite{Lam}. To us this means that every finitely generated and projective module over $\A$ is free. Later this was generalized to an even bigger class of rings by Gubaladze in \cite{Gubeladze}, see also \cite{Lam2}. However, we shall now see that the theorem by Swan has consequences for the structure of finitely generated shift-invariant subspaces in $L^2(\R^n)$.

\begin{theorem}\label{main}
If  
$\zeta_1, \dots, \zeta_d$ and $ \tilde \zeta_1, \dots, \tilde \zeta_d$ ore compactly supported and form dual $G$-frames for a closed subspace $\mathcal K \subset L^2(\R^n)$, then 
%
there exist bi-orthogonal $G$-frames 
 $\eta_1, \dots, \eta_r$ and $ \tilde \eta_1, \dots, \tilde \eta_d$ with compact supports in $\K$. 
\begin{proof}
  Let $E \subset \Lambda$ denote the $\A$-module generated by $\zeta_1, \dots , \zeta_r$. 
  By Theorem \ref{frametheorem}, $E$ is self dual and projective. However, since $E$ is a finitely generated projective module over $\A$, it is necessarily free. If we apply Lemma \ref{moduleframelemma}, we see that since $E$ is self dual, there exist $\eta_1, \dots, \eta_r$ and $ \tilde \eta_1, \dots, \tilde \eta_r$ in $E$ such that 
  $\langle \eta_i, \tilde \eta_j \rangle_\A = \delta_{i,j}$. If we apply Corollary \ref{biortcor}, we see that 
  $\eta_1, \dots, \eta_r$ and $ \tilde \eta_1, \dots, \tilde \eta_r$ form bi-orthogonal $G$-frames for $\K$.
 %
%

\end{proof}
\end{theorem}

\section{Two  applications to multiresolution analysis.} 
Let $G = \Z^n$, let $\nu$ denote the group action defined in the previous section and let $\Lambda$ denote set of compactly supported elements in $L^2(\R^n)$.  Moreover, let $A \in GL_n( \R)$, such that $AG \subset G$. We assume that 
every eigenvalue of $A$ has absolute value strictly greater than $1$. The determinant of $A$
is an integer. Let $q = | \det A |$. We define a unitary operator on $L^2( \R^n)$ 
with the equation: $U \zeta (t) = \sqrt q \zeta(At)$. 
\begin{definition} \label{MRAdef}
We will say that $U$ and $\phi_1, \dots, \phi_d, \tilde \phi_1, \dots, \tilde \phi_d \in L^2(\R^n)$ generate a \textit{multiresolution analysis}, MRA for short, if the following properties are satisfied:
\begin{enumerate}
\item 
	The families $\phi_1, \dots, \phi_d$ and $\tilde \phi_1, \dots, \tilde \phi_d$ form dual $G$-frames for 
	a closed subspace $\V \subset L^2(\R^n)$.
\item The subspace $U \V$ contains $\V$.
\item The union $\cup_{k \in \Z} U^k \V$ is dense in $L^2(\R^n)$.
\item The intersection $\cap_{k \in \Z} U^k \V$ equals the trivial subspace $\{0\}$.
\end{enumerate} 
The subspace $\V$ is called the \textit{scaling space} and the subspace 
$U \V \ominus \V$ is called the \textit{wavelet space}.
The families $\phi_1, \dots, \phi_d$ and $\tilde \phi_1, \dots, \tilde \phi_d$ are said to be \textit{dual scaling families} for the MRA. 
Two families $\psi_1, \dots , \psi_r \in U\V$ and $\tilde \psi_1, \dots , \tilde \psi_r \in U\V$ are said to be \textit{dual MRA wavelet families} if they form dual $G$-frames for the wavelet space.  If the associated frames are bi-orthogonal, we will say that they form \textit{bi-orthogonal MRA wavelet families}. 
\end{definition}
\subsection{Bi-orthogonal wavelet families.}
Suppose  
	$\phi_1, \dots, \phi_d$ and $\tilde \phi_1, \dots, \tilde \phi_d$ are compactly supported and that they 
	define dual scaling families for an MRA with scaling space $\V$. 
	Now $E = \Lambda \cap \V$ coincides with the compactly supported elements in $\V$. Theorem \ref{main} tells us that there exist bi-orthogonal $G$-frames for $\V$ in $\Lambda$. The next theorem tells us in addition that $U\V \cap \Lambda$  contains bi-orthogonal $G$-frames for the wavelet space.  

\begin{theorem}
	Suppose  
	$\phi_1, \dots, \phi_d$ and $\tilde \phi_1, \dots, \tilde \phi_d$ are compactly supported and that they 
	define  dual scaling families for an MRA with scaling space $\V$ as in Definition \ref{MRAdef}.
	
	The  wavelet space contains compactly supported bi-orthogonal MRA wavelet families.
  \begin{proof}
  		 Let $E$ denote the space of compactly supported elements in $\V$, i.e. $E$ is the $\A$-submodule in $\Lambda$ generated by $\phi_1, \dots, \phi_d$. The intersection $U \V \cap \Lambda = UE$ defines an $\A$-submodule in $\Lambda$ that contains $E$ as an $\A$-submodule.
  		 
  	 By Theorem \ref{main}, we can choose  
  	 $\phi_1, \dots, \phi_d$ and $\tilde \phi_1, \dots, \tilde \phi_d$
  	 to be bi-orthogonal and dual $G$-frames for $\V$. Now 
  	 $\phi_1, \dots, \phi_d$ and $\tilde \phi_1, \dots, \tilde \phi_d$ form bi-orthogonal $\A$-module frames 
  	 for $E$.

  	 Let $g_1, \dots, g_q$ be a system of representatives of the cosets in $G/AG$. Since $\nu_g U = U\nu_{Ag}$ for every $g \in G$, we see immediately that 
  	 $\{ U\nu_{g_i} \phi_j | 1 \leq i \leq q, 1 \leq j  \leq d\}$
  and  $\{ U\nu_{g_i} \tilde \phi_j | 1 \leq i \leq q, 1 \leq j  \leq d\}$ form bi-orthogonal $G$-frames for 
  the subspace $U\V$. By  Corollary \ref{biortcor} they also form bi-orthogonal $\A$-module frames for the free $\A$-module $UE$.  
  Let $P\in B( L^2( \R^n)$) denote the orthogonal projection onto $\V$.	  
  We see that 
  \begin{align*}
  	P \zeta = \sum_{i = 1}^d \langle \zeta , \tilde \phi_i \rangle_\A \phi_i = 
  	\sum_{i = 1}^d \langle \zeta , \phi_i \rangle_\A \tilde \phi_i
  \end{align*}
  for every $\zeta \in \Lambda$. The map $P$ is $\A$-linear, self adjoint and idempotent on $\Lambda$ such that 
  $P \Lambda = E$. Now 
  $\{ PU\nu_{g_i} \phi_j | 1 \leq i \leq q, 1 \leq j  \leq d\}$
  and  $\{ PU\nu_{g_i} \tilde \phi_j | 1 \leq i \leq q, 1 \leq j  \leq d\}$
  form dual $\A$-module frames for the module $(1-P)UE$.
  By Lemma \ref{moduleframelemma} we see that the submodule $(1-P)UE\subset\Lambda$ is finitely generated, self dual and projective. 
  Finally,  Theorem \ref{main} tells us that there exist compactly supported bi-orthogonal $G$-frames 
  $\psi_1,\dots,\psi_r$ and $\tilde\psi_1,\dots,\tilde\psi_r$ in  $\overline{(1-P)UE}=U\V\ominus\V$.

  \end{proof}
\end{theorem}

\subsection{Symmetric MRA wavelets.}
There are several examples of multiresolution analyses in $L^2(\R)$ generated  
by dual and compactly supported bi-orthogonal $G$-frames $\phi$ and $\tilde  \phi $  such that $\phi(x) = \phi(-x) \in \R$ and 
$\tilde \phi(x) = \tilde \phi(-x) \in \R$ for almost every $x$. Moreover, these multiresolution analyses allows 
bi-orthogonal and compactly supported wavelets $\psi$ and $\tilde \psi$ such that 
$\psi(x + 2\-  ) = \psi(-x+ 2\- ) \in \R$ and 
$\tilde \psi(x+ 2\- ) = \tilde \psi(-x +2\- ) \in \R$ for almost every $x$, see \cite{Daubechies10Lectures} and  \cite{Daubechies2}.

We will be interested in how one can construct bi-orthogonal and compactly supported higher dimensional MRA wavelets with symmetries if we have compactly supported scaling functions with some symmetry properties. 
The following definition will turn out to be important: 
\begin{definition}
We will say that a finite subgroup $H \subset GL_n(\Z)$  \textit{is affiliated to the dilation} $A$ if it satisfies the following properties:  
\begin{enumerate}
\item $h A = A h$, \label{comm}
\item $( h-I) G \subset AG$, \label{triv}
\end{enumerate}
for every $h \in H$.
\end{definition}
Note that if $n = 1$ and $H$ is non trivial, then  $A = \pm 2$ and $H = \{\pm 1\}$ are the only examples. We classified in \cite{roysland1} 
every dilation with a nontrivial affiliated group up to similarity when $n =2$.  
This classification shows that if that if $n = 2$ and $q = 2$ then either $H \simeq \Z/ 2\Z$ or $H \simeq \Z/ 4 \Z$. 
If $H \simeq \Z^/ 2 \Z$, then 
$H$ is generated by the matrix $-I_2$
and $A$ is similar to one of the following matrices:
$$
 \begin{pmatrix}
   0 &  2 \\
  \pm 1 &  0 \\
 \end{pmatrix},
 \pm \begin{pmatrix}
   0  &  2 \\
  -1 &  1 \\
 \end{pmatrix},
 \pm\begin{pmatrix}
   1 &  -1 \\
  1 &  1 \\
 \end{pmatrix}.
$$
Moreover, if $H \simeq \Z/ 4 \Z$ and $h$ generates $H$, then there exists an $S \in GL_2(\Z)$ such that 
\begin{align*}
ShS^{-1} =  
\begin{pmatrix}
   0 & 1 \\ -1&  0 \\
\end{pmatrix} 
 ~\text{ and } ~ 
 SAS^{-1} = \pm\begin{pmatrix}
   1 &  -1 \\
  1 &  1 \\
 \end{pmatrix}.
\end{align*}

\begin{theorem}
	Suppose $q =2$ and suppose $\phi$ and $\tilde \phi$ have compact supports and generate an MRA with dilation $A$ and 
	scaling space $\V$. 
	If $g_1 \in G \setminus AG$, then the vectors
	\begin{align*}
	\psi & = \langle U \nu_{g_1} \phi, \tilde \phi \rangle_\A U \phi - \langle U \phi , \tilde \phi \rangle_\A U \nu_{g_1} \phi \\
	\tilde \psi & =
	\langle U \nu_{g_1} \tilde \phi,  \phi \rangle_\A U \tilde \phi - \langle U \tilde \phi ,  \phi \rangle_\A U \nu_{g_1} \tilde \phi \\
	\end{align*}
	define bi-orthogonal MRA wavelets in $\Lambda \cap U \V$.  
	Moreover, if $H$ is affiliated to $A$ and 
	  
	  $\phi (hx)  = \phi(x)$ and $\tilde \phi (hx)    = \tilde \phi(x)$
	  for every 
	$h \in H$ and a.e. $x \in \R^n$, then 
	\begin{align} \label{symwav}
	\psi ( h x + A^{-1} g_1)  = \psi (  x + A^{-1} g_1)~ \text{ and } ~
	\tilde \psi ( h x + A^{-1} g_1)  = \tilde \psi (  x + A^{-1} g_1)
	\end{align}
	for every $h \in H$ and a.e. $x \in \R^n$.
	\begin{proof}
		First, let $E$ denote the $\A$-module $\Lambda \cap \V$. Note that $U \zeta \in \Lambda$ if and only if $U \zeta \in \Lambda$. 
		This and the assumption $\V \subset U \V$ implies that  $E \subset U E \subset \Lambda$ and $U \V \cap \Lambda  = U E$.
		We have that $\{ U \nu_g \phi\}_{g \in G}$ and 	$\{ U \nu_g \tilde \phi\}_{g \in G}$ 	form bi-orthogonal frames for 
		the space $U\V$. The element $g_1$ represents the only non-trivial coset in $G$. Since $\nu_gU = \nu_{Ag}$ for every $g \in G$, we see that  
		$\{\nu_g U \nu_f \phi | f \in \{  0,g_1\} , g \in G  \}$ and $\{\nu_g U \nu_f \tilde \phi | f \in \{  0,g_1\} , g \in G  \}$ 
		form bi-orthogonal frames for $U \V$. Now $\{U  \phi, U \nu_{g_1} \phi\}$ and $\{U  \tilde \phi, U \nu_{g_1} \tilde \phi\}$ 
		form dual module frames for the $\A$-module $UE$. To simplify our notation slightly, we let $\zeta_1 = U \phi$, 
		$\zeta_2 = U \nu_{g_1} \phi$, $\tilde \zeta_1 = U \tilde \phi$ and $\tilde \zeta_2 = U \nu_{g_1} \tilde \phi$.
		
		We will  first check that $\{\phi, \psi\}$ and $\{ \tilde \phi , \tilde \psi\}$ form  dual module frames   for $UE$:
		Since 	$\{\zeta_1, \zeta_2\}$ and $\{\tilde \zeta_1, \tilde \zeta_2\}$  form dual module bases for 
		$UE$, we obtain the following identities:
		\begin{align*}
			\langle \tilde \zeta_1, \psi \rangle_\A \langle \tilde \psi, \zeta_1 \rangle_\A & =
			\langle \tilde \phi, \zeta_2\rangle_\A \langle \tilde \zeta_2, \phi \rangle_\A \\
			\langle \tilde \zeta_1, \psi\rangle_\A \langle \tilde \psi, \zeta_2 \rangle_\A & = - \langle \tilde \zeta_1, \phi \rangle_\A 
			\langle \tilde \phi, \zeta_2 \rangle_\A \\
			\langle \tilde \zeta_2, \psi \rangle_\A \langle \tilde \psi, \zeta_1 \rangle_\A & = - \langle \tilde \zeta_2 , \phi \rangle_\A
			\langle \tilde \phi,  \zeta_1 \rangle_\A \\
			\langle \zeta_2, \psi \rangle_\A \langle \tilde \psi, \zeta_2 \rangle_\A & =
			\langle \tilde \zeta_1, \phi \rangle_A \langle \tilde \phi, \zeta_1 \rangle_\A.  
		\end{align*}
		Putting these together, we see that
		\begin{align} \label{biorteq}
			\langle \tilde \zeta_i, \phi \rangle_\A \langle \tilde \phi , \zeta_j \rangle_\A + 
			\langle \tilde \zeta_i, \psi \rangle_\A \langle \tilde \psi , \zeta_j \rangle_\A  = 
			\delta_{i,j}. 
		\end{align}
		If we let $\zeta, \eta \in UE$ and apply \eqref{biorteq},  then  we obtain 
		\begin{align*}
			& \langle \zeta, \phi \rangle_\A  \langle \tilde \phi ,  \eta \rangle_\A + 
			\langle \zeta, \psi \rangle_\A  \langle \tilde \psi ,  \eta \rangle_\A
			 \\= &  \sum_{i,j} \langle \zeta,  \zeta_i \rangle_\A \langle \tilde \zeta_i, \phi \rangle_\A
			 \langle \tilde \phi, \zeta_j \rangle_\A \langle \tilde  \zeta_j ,\eta  \rangle_\A   
			 + 
			 \sum_{i,j} \langle \zeta,  \zeta_i \rangle_\A \langle \tilde \zeta_i, \psi \rangle_\A
			 \langle \tilde \psi, \zeta_j \rangle_\A \langle \tilde  \zeta_j ,\eta  \rangle_\A   \\
			  = &  \sum_{i,j} \langle \zeta,  \zeta_i \rangle_\A \big( \langle \tilde \zeta_i, \phi \rangle_\A
			 \langle \tilde \phi, \zeta_j \rangle_\A + 
			 \langle \tilde \zeta_i, \psi \rangle_\A
			 \langle \tilde \psi, \zeta_j \rangle_\A \big) \langle \tilde  \zeta_j ,\eta  \rangle_\A \\
			  = & \sum_i \langle \zeta,  \zeta_i \rangle_\A    \langle \tilde  \zeta_i ,\eta  \rangle_\A  = \langle \zeta, \eta \rangle_\A,
		\end{align*}
	so 	$\{\phi, \psi\}$ and $\{ \tilde \phi , \tilde \psi\}$ form  dual module frames   for $UE$.
	
	Next, we compute that 
	\begin{align*}
		\langle \phi, \tilde \psi \rangle_\A = \langle \phi, \tilde \zeta_1 \rangle_\A \langle \phi, \tilde \zeta_2\rangle_\A 
		- \langle \phi, \tilde \zeta_2\rangle_\A \langle \phi , \tilde \zeta_1 \rangle_\A = 0.
	\end{align*}
	An equivalent computation shows that  $\langle \tilde \phi,  \psi \rangle_\A = 0$. 
	Finally, we compute that 
	\begin{align*}
		\langle \psi, \tilde \psi \rangle_\A & = \langle \langle \zeta_2 , \tilde \phi \rangle_\A \zeta_1 - \langle \zeta_1, \tilde \phi \rangle_\A \zeta_2, \langle \tilde \zeta_2 ,  \phi \rangle_\A \tilde \zeta_1 - \langle \tilde \zeta_1,  \phi \rangle_\A \tilde \zeta_2 \rangle_\A \\
		& = \sum_i \langle \phi, \tilde \zeta_i \rangle_\A \langle \zeta_i , \tilde \phi \rangle_\A = \langle \phi, \tilde \phi \rangle_\A = 1, 
	\end{align*}
	so $\{\phi, \psi\}$ and $\{ \tilde \phi , \tilde \psi\}$ form  dual module bases   for $UE$.
	By Corollary \ref{biortcor}, we see that $\psi$ and $\tilde \psi$ form bi-orthogonal and compactly supported MRA wavelets in $UE$. 
	For every $h \in H$, we define a unitary operator on $L^2( \R^n)$ with  the following equation:
	$
		W_h \xi (x) = \xi( h^{-1} x). 
	$
	These unitaries define a unitary representation of $H$ on $L^2 ( \R^n )$. Whenever $\eta_1, \eta_2 , \zeta \in \Lambda$, then
	\begin{align*}
		\langle W_h \eta_1, W_h \eta_2 \rangle_\A \zeta  & = \sum_{g \in G} \langle W_h \eta_1, \nu_g W_h \eta_2 \rangle \nu_g \zeta 
	= \sum_{g \in G} \langle \eta_1, W_h ^{-1} \nu_g W_h \eta_2 \rangle \nu_g \zeta \\
	& = \sum_{g \in G} \langle \eta_1, \nu_g \eta_2 \rangle W_h \nu_g W_h^{-1} \zeta 
		 = W_h \langle \eta_1, \eta_2 \rangle_\A W_h^{-1} \zeta.
	\end{align*}
	We have that $W_h U = U W_h$ for every $h \in H$, since $h A = Ah$ for every $h \in H$. 
	Moreover, since $(h -1)G \subset AG$ for every $h \in H$,  there exists a $g_h \in G$ for every $h \in H$ s.t.
	$$U W_h \nu_{g_1} W_h ^{-1} \nu_{g_1}^{-1} U^{-1} = \nu_{g_h}.$$
	We also see that $U\nu_{g_1}  U^{-1}$ commutes with every element in $\pi_0( \A)$. 	
	Now
	\begin{align*}
		W_h U \nu_{g_1}^{-1}  U^{-1} \psi   = &  
		W_h U \nu_{g_1}^{-1}  U^{-1} \langle U \nu_{g_1} \phi, \tilde \phi \rangle_\A U \phi - 
		W_h U \nu_g^{-1}  U^{-1} \langle U \phi , \tilde \phi \rangle_\A U \nu_{g_1} \phi \\
		 = &  W_h U \nu_{g_1}^{-1}  U^{-1} W_h^{-1} W_h \langle U \nu_{g_1}  \phi, \tilde \phi \rangle_\A W_h^{-1} U \phi 
		- W_h \langle U \phi, \tilde \phi \rangle_\A U \phi \\
		 = & W_h U \nu_{g_1}^{-1}  U^{-1} W_h^{-1}  \langle U W_h \nu_{g_1}  W_h^{-1}  \phi, \tilde \phi \rangle_\A  U \phi 
		- W_h \langle U \phi, \tilde \phi \rangle_\A W_h^{-1} U \phi \\
		 = & W_h U \nu_{g_1}^{-1}  U^{-1} W_h^{-1}  \langle \nu_{g_h}  U \nu_{g_1} \phi, \tilde \phi \rangle_\A  U \phi 
		- \langle W_h U \phi, W_h \tilde \phi \rangle_\A  U \phi \\
		 = & W_h U \nu_{g_1}^{-1}  U^{-1} W_h^{-1} \nu_{g_h}  \langle   U \nu_{g_1} \phi, \tilde \phi \rangle_\A  U \phi 
		- \langle U \phi,  \tilde \phi \rangle_\A  U \phi \\
		 = &  U \nu_{g_1}^{-1} U^{-1}  \langle   U \nu_{g_1} \phi, \tilde \phi \rangle_\A  U \phi 
		- \langle U \phi,  \tilde \phi \rangle_\A  U \phi \\
		 = & U \nu_{g_1}^{-1}  U^{-1}  (   \langle U \nu_{g_1} \phi, \tilde \phi \rangle_\A U \phi - \langle U \phi , 
		\tilde \phi \rangle_\A U \nu_{g_1} \phi                        )
		 = U \nu_{g_1}^{-1}  U^{-1} \psi, 
	\end{align*}
	i.e. $\psi ( h x + A^{-1} g_1)  = \psi (  x + A^{-1} g_1)$ for almost every $x \in \R^n$. 
	The analogous statement about $\tilde \psi$ follows from an identical computation. 
	\end{proof}
\end{theorem}

\begin{acknowledgements}The author is  pleased to acknowledge helpful
discussions with Ola Bratteli and Martin Gulbrandsen.
\end{acknowledgements}

\bibliographystyle{alpha}
\bibliography{invariant}

\begin{thebibliography}{CDF92}

\bibitem[BR87]{BratteliRobinson}
Ola Bratteli and Derek~W. Robinson.
\newblock {\em Operator algebras and quantum statistical mechanics. 1}.
\newblock Texts and Monographs in Physics. Springer-Verlag, New York, second
  edition, 1987.
\newblock $C\sp \ast$- and $W\sp \ast$-algebras, symmetry groups, decomposition
  of states.

\bibitem[CDF92]{Daubechies2}
A.~Cohen, Ingrid Daubechies, and J.-C. Feauveau.
\newblock Biorthogonal bases of compactly supported wavelets.
\newblock {\em Comm. Pure Appl. Math.}, 45(5):485--560, 1992.

\bibitem[Dau92]{Daubechies10Lectures}
Ingrid Daubechies.
\newblock {\em Ten lectures on wavelets}, volume~61 of {\em CBMS-NSF Regional
  Conference Series in Applied Mathematics}.
\newblock Society for Industrial and Applied Mathematics (SIAM), Philadelphia,
  PA, 1992.

\bibitem[Dav96]{byexamples}
Kenneth~R. Davidson.
\newblock {\em {$C\sp *$}-algebras by example}, volume~6 of {\em Fields
  Institute Monographs}.
\newblock American Mathematical Society, Providence, RI, 1996.

\bibitem[DR07]{DutRoy1}
Dorin~Ervin Dutkay and Kjetil R{\o}ysland.
\newblock The algebra of harmonic functions for a matrix-valued transfer
  operator.
\newblock {\em J. Funct. Anal.}, 252(2):734--762, 2007.

\bibitem[FL02]{franklarson}
Michael Frank and David~R. Larson.
\newblock Frames in {H}ilbert {$C\sp \ast$}-modules and {$C\sp \ast$}-algebras.
\newblock {\em J. Operator Theory}, 48(2):273--314, 2002.

\bibitem[Gub88]{Gubeladze}
I.~Dzh. Gubeladze.
\newblock The {A}nderson conjecture and a maximal class of monoids over which
  projective modules are free.
\newblock {\em Mat. Sb. (N.S.)}, 135(177)(2):169--185, 271, 1988.

\bibitem[HL00]{hanlarson}
Deguang Han and David~R. Larson.
\newblock Frames, bases and group representations.
\newblock {\em Mem. Amer. Math. Soc.}, 147(697):x+94, 2000.

\bibitem[Lam78]{Lam}
T.~Y. Lam.
\newblock {\em Serre's conjecture}.
\newblock Springer-Verlag, Berlin, 1978.
\newblock Lecture Notes in Mathematics, Vol. 635.

\bibitem[Lam06]{Lam2}
T.~Y. Lam.
\newblock {\em Serre's problem on projective modules}.
\newblock Springer Monographs in Mathematics. Springer-Verlag, Berlin, 2006.

\bibitem[Lan95]{Lance}
E.~C. Lance.
\newblock {\em Hilbert {$C\sp *$}-modules}, volume 210 of {\em London
  Mathematical Society Lecture Note Series}.
\newblock Cambridge University Press, Cambridge, 1995.
\newblock A toolkit for operator algebraists.

\bibitem[Lan02]{Lang}
Serge Lang.
\newblock {\em Algebra}, volume 211 of {\em Graduate Texts in Mathematics}.
\newblock Springer-Verlag, New York, third edition, 2002.

\bibitem[LR07]{Larsen}
Nadia~S. Larsen and Iain Raeburn.
\newblock Projective multi-resolution analyses arising from direct limits of
  {H}ilbert modules.
\newblock {\em Math. Scand.}, 100(2):317--360, 2007.

\bibitem[Pac07]{Packer}
Judith~A. Packer.
\newblock Projective multiresolution analyses for dilations in higher
  dimensions.
\newblock {\em J. Operator Theory}, 57(1):147--172, 2007.

\bibitem[Ped89]{Pedersen}
Gert~K. Pedersen.
\newblock {\em Analysis now}, volume 118 of {\em Graduate Texts in
  Mathematics}.
\newblock Springer-Verlag, New York, 1989.

\bibitem[PR03]{Pr2}
Judith~A. Packer and Marc~A. Rieffel.
\newblock Wavelet filter functions, the matrix completion problem, and
  projective modules over {$C(\Bbb T\sp n)$}.
\newblock {\em J. Fourier Anal. Appl.}, 9(2):101--116, 2003.

\bibitem[PR04]{PR1}
Judith~A. Packer and Marc~A. Rieffel.
\newblock Projective multi-resolution analyses for {$L\sp 2({\Bbb R}\sp 2)$}.
\newblock {\em J. Fourier Anal. Appl.}, 10(5):439--464, 2004.

\bibitem[R{\o}y08]{roysland1}
Kjetil R{\o}ysland.
\newblock Symmetries in projective multiresolution analyses.
\newblock {\em J. Fourier Anal. Appl.}, 14(2):267--285, 2008.

\bibitem[Swa78]{Swan}
Richard~G. Swan.
\newblock Projective modules over {L}aurent polynomial rings.
\newblock {\em Trans. Amer. Math. Soc.}, 237:111--120, 1978.

\bibitem[Woo04]{wood}
Peter~John Wood.
\newblock Wavelets and {H}ilbert modules.
\newblock {\em J. Fourier Anal. Appl.}, 10(6):573--598, 2004.

\end{thebibliography}

\end{document}